\input amsppt.sty
\magnification=\magstep1
\hsize=30truecc
\vsize=22.2truecm
\baselineskip=16truept
\nologo
\pageno=1
\TagsOnRight
\topmatter
\def\Z{\Bbb Z}
\def\N{\Bbb N}

\def\C{\Bbb C}

\def\l{\left}
\def\r{\right}
\def\bg{\bigg}
\def\({\bg(}
\def\[{\bg[}
\def\){\bg)}
\def\]{\bg]}
\def\t{\text}
\def\f{\frac}

\def\bi{\binom}

\def\bi{\binom}

\def\Proof{\noindent{\it Proof}}
\def\Remark{\medskip\noindent{\it  Remark}}

\topmatter \hbox{Last modified: August 20, 2004; {\tt
arXiv:math.NT/0408223}.}
\bigskip
\title {New identities involving Bernoulli and Euler polynomials. II}\endtitle
\rightheadtext{Identities involving Bernoulli and Euler
polynomials. II}
\author {Zhi-Wei Sun and Hao Pan}\endauthor
\affil Department of Mathematics, Nanjing University
\\Nanjing 210093, The People's Republic of China
\\ {\tt zwsun\@nju.edu.cn}\ \ \quad\tt{haopan79\@yahoo.com.cn}
\endaffil
\abstract We derive several symmetric identities for
Bernoulli and Euler polynomials which imply some known identities.
Our proofs depend on the new technique developed in part I and
some identities obtained in [European J. Combin. 24(2003), 709--718].
\endabstract
\thanks  2000 {\it Mathematics Subject Classification}.
Primary 11B68; Secondary 05A19.
\newline\indent The first author is responsible for all the communications
and supported by
the National Natural Science Foundation of P. R. China.
\endthanks
\endtopmatter
\document
\hsize=30truecc
\vsize=22.2truecm
\baselineskip=16truept

\heading 1. Introduction\endheading

Bernoulli polynomials $B_n(x)\ (n=0,1,2,\ldots)$ and Euler polynomials
$E_n(x)\ (n=0,1,2,\ldots)$ are defined by power series
$$\frac{ze^{xz}}{e^z-1}=\sum_{n=0}^{\infty}B_n(x)\frac{z^n}{n!}
\ \ \ \t{and}\ \ \ \frac{2e^{xz}}{e^z+1}=\sum_{n=0}^{\infty}E_n(x)\frac{z^n}{n!}.$$
Those $B_n=B_n(0)$ and $E_n=2^nE_n(1/2)$ are called Bernoulli numbers and Euler numbers
respectively.
Here are some well-known properties of Bernoulli and Euler polynomials (cf. [S1]):
$$\gather B_n(1-x)=(-1)^nB_n(x),\ B_n(x+y)=\sum_{k=0}^n\bi{n}{k}B_{k}(x)y^{n-k},
\\E_n(1-x)=(-1)^nE_n(x),\ E_n(x+y)=\sum_{k=0}^n\bi{n}{k}E_{k}(x)y^{n-k}.\endgather$$
Also,
$$\Delta(B_n(x))=nx^{n-1}\ \ \t{and}\ \ \Delta^*(E_n(x))=2x^n,$$
where
$\Delta(f(x))=f(x+1)-f(x)$ and $\Delta^*(f(x))=f(x+1)+f(x)$.
($\Delta$ is usually called the difference operator.)
For $n\in \Z^+=\{1,2,3,\ldots\}$, we also have
$B_n'(x)=nB_{n-1}(x)$ and $E_n'(x)=nE_{n-1}(x)$.

In 1979 C. F. Woodcock [W] discovered that $A_{m-1,\,n}=A_{n-1,\,m}$ for $m,n\in\Z^+$
where
$$A_{m,\,n}=\f1n\sum_{k=1}^n\bi nk(-1)^kB_{m+k}B_{n-k}.$$

In 2003 the first author [S2] obtained the following result concerning Bernoulli and Euler polynomials.

\proclaim{Theorem 0 {\rm ([S2, Theorem 1.2])}} Let $m,n\in\N=\{0,1,2,\ldots\}$ and $x+y+z=1$.
Then
$$\aligned&(-1)^m\sum_{k=0}^m\bi mkx^{m-k}\f{B_{n+k+1}(y)}{n+k+1}
+(-1)^n\sum_{k=0}^n\bi nkx^{n-k}\f{B_{m+k+1}(z)}{m+k+1}
\\&\qquad\qquad=\f{(-x)^{m+n+1}}{(m+n+1)\bi {m+n}n}
\endaligned$$
and
$$(-1)^m\sum_{k=0}^m\bi mkx^{m-k}B_{n+k}(y)=(-1)^n\sum_{k=0}^n\bi nkx^{n-k}B_{m+k}(z).$$
Also, we can replace all the Bernoulli polynomials in the last two identities by corresponding
Euler polynomials.
\endproclaim

In part I ([PS]) we obtained  polynomial versions of
Miki's and Matiyasevich's curious identities for Bernoulli numbers,
the new method developed there involves differences and derivatives of polynomials.

In this paper, by using the technique in part I, we will
prove some identities (for Bernoulli and Euler polynomials)
related to both Woodcock's result and Theorem 0.

Here is our central result in this paper.

\proclaim{Theorem 1} Let $m,n\in\N$ and $x+y+z=1$. Then
$$\aligned
&(-1)^m\sum_{k=0}^{m}\bi mk\frac{B_{m-k+1}(x)}{m-k+1}\cdot\f{B_{n+k+1}(y)}{n+k+1}
\\&+(-1)^n\sum_{k=0}^{n}\bi nk\frac{B_{n-k+1}(x)}{n-k+1}\cdot\f{B_{m+k+1}(z)}{m+k+1}
\\=&\f{(-1)^{m+n+1}}{(m+n+1)\bi {m+n}n}\cdot\f{B_{m+n+2}(x)}{m+n+2}-\f{B_{m+1}(z)}{m+1}\cdot\f{B_{n+1}(y)}{n+1}
\\&+\f{(-1)^{m+1}}{m+1}\cdot\f{B_{m+n+2}(y)}{m+n+2}+\f{(-1)^{n+1}}{n+1}\cdot\f{B_{m+n+2}(z)}{m+n+2}.
\endaligned\tag 1$$
Also,
$$\aligned
&(-1)^m\sum_{k=0}^m\bi mkE_{m-k}(x)\f{B_{n+k+1}(y)}{n+k+1}
\\&+(-1)^n\sum_{k=0}^n\bi nkE_{n-k}(x)\f{B_{m+k+1}(z)}{m+k+1}
\\=&\f{(-1)^{m+n+1}E_{m+n+1}(x)}{(m+n+1)\bi{m+n}n}-\f{E_m(z)E_n(y)}2
\endaligned\tag 2$$
and
$$\aligned&\f{(-1)^m}2\sum_{k=0}^m\bi mkE_{m-k}(x)\f{E_{n+k+1}(y)}{n+k+1}
\\&-(-1)^n\sum_{k=0}^n\bi nk\f{B_{n-k+1}(x)}{n-k+1}\cdot\f{E_{m+k+1}(z)}{m+k+1}
\\=&\f{(-1)^{m+n}}{(m+n+1)\bi{m+n}n}\cdot\f{B_{m+n+2}(x)}{m+n+2}
+\f{(-1)^n}{n+1}\cdot\f{E_{m+n+2}(z)}{m+n+2}
\\&-\f 1{n+1}\sum_{k=0}^m\f{\bi mk}{\bi {n+k+1}k}E_{m-k}(z)\f{B_{n+k+2}(y)}{n+k+2}.
\endaligned\tag3$$
\endproclaim

\Remark\ 1. Fix $y$ and replace $z$ in (1) by $1-x-y$.
Then, by taking differences of both sides of (1) with respect to $x$,
we can get the first identity in Theorem 0. Similarly, other identities in Theorem 0
are also implied by Theorem 1.
\medskip

Clearly (1) has the following equivalent form:
$$\align&\f{(-1)^m}{m+1}\sum_{k=0}^{m+1}\bi {m+1}kB_{m+1-k}(x)\f{B_{n+1+k}(y)}{n+1+k}
\\&+\f{(-1)^n}{n+1}\sum_{k=0}^{n+1}\bi {n+1}kB_{n+1-k}(x)\f{B_{m+1+k}(z)}{m+1+k}
\\=&\f{(-1)^{m+n+1}m!n!}{(m+n+2)!}B_{m+n+2}(x)-\f{B_{m+1}(z)}{m+1}\cdot\f{B_{n+1}(y)}{n+1}.
\endalign$$
So, if $m,n\in\Z^+$ and $x+y+z=1$, then we have
$$\aligned &\f{(-1)^m}m\sum_{k=0}^m\bi mkB_{m-k}(x)\f{B_{n+k}(y)}{n+k}
+\f{(-1)^n}n\sum_{k=0}^n\bi nkB_{n-k}(x)\f{B_{m+k}(z)}{m+k}
\\&\qquad=\f{(-1)^{m+n}(m-1)!(n-1)!}{(m+n)!}B_{m+n}(x)+\f{B_m(z)}m\cdot\f{B_n(y)}n.
\endaligned\tag1$'$ $$

\proclaim{Corollary 1} Let $x+y+z=1$. Given $m,n\in\Z^+$ we have the following identities:
$$\aligned&\f{(-1)^m}m\sum_{k=0}^{m}\bi mkB_{m-k}(x)B_{n-1+k}(y)-\f{B_m(z)}mB_{n-1}(y)
\\=&\f{(-1)^n}n\sum_{k=0}^{n}\bi nkB_{n-k}(x)B_{m-1+k}(z)-\f{B_n(y)}nB_{m-1}(z),
\endaligned\tag4$$
$$\aligned &(-1)^m\sum_{k=0}^m\bi mkE_{m-k}(x)B_{n+k}(y)-\f m2E_{m-1}(z)E_n(y)
\\&=(-1)^n\sum_{k=0}^n\bi nkE_{n-k}(x)B_{m+k}(z)-\f n2E_{n-1}(z)E_m(y)
\endaligned\tag5$$
and
$$\aligned
&\f{(-1)^m}2\sum_{k=0}^m\bi mkE_{m-k}(x)E_{n-1+k}(y)
\\=&\f{(-1)^n}n\sum_{k=0}^n\bi nkB_{n-k}(x)E_{m+k}(z)-\f{B_n(y)}nE_m(z).
\endaligned\tag 6$$
\endproclaim
\Proof. Replacing $z$ in ($1'$) by $1-x-y$ and taking partial derivatives with respect to $y$,
we then obtain (4) from ($1'$).
Identity (5) can be easily deduced from (2) by
taking partial derivatives with respect to $y$.
Similarly, (6) follows from (3) with $n$ replaced by $n-1$.
 \qed

\proclaim{Corollary 2} For  $m,n\in\Z^+$ we have
$$A_{m-1,\,n}(t)=A_{n-1,\,m}(t)\ \ \t{and}\ \ C_{m,\,n}(t)=C_{n,\,m}(t),\tag7$$
where
$$A_{m,\,n}(t)=\f1n\sum_{k=0}^n\bi nk(-1)^kB_{m+k}(t)B_{n-k}(2t)-B_m(t)\f{B_n(t)}n\tag8$$
and
$$C_{m,\,n}(t)=\sum_{k=0}^n\bi nk(-1)^kB_{m+k}(t)E_{n-k}(2t)-\f n2E_m(t)E_{n-1}(t).\tag9$$
\endproclaim
\Proof. Just apply (4) and (5) with $x=1-2t$ and $y=z=t$. \qed

\Remark\ 2. The first equality in (7) is an extension of Woodcock's identity.

\proclaim{Corollary 3} If  $m,n\in\Z^+$ then
$$\aligned&\f1m\sum_{k=0}^m\bi mk(-1)^k\f{(1-2^{n+k})B_{n+k}}{n+k}(1-2^{m-k})B_{m-k}
\\&\ =\f1n\sum_{k=1}^n\bi nk(-1)^k\f{(1-2^{m+k})B_{m+k}}{m+k}B_{n-k}.
\endaligned\tag10$$
\endproclaim
\Proof. For $l\in\N$ it is well known (cf. [S1]) that
$$E_l(x)=\f 2{l+1}\(B_{l+1}(x)-2^{l+1}B_{l+1}\l(\f x2\r)\),$$
thus $(-1)^lE_l(1)=E_l(0)=2(1-2^{l+1})B_{l+1}/(l+1)$.
Applying (6) with  $x=1$ and $y=z=0$ and replacing $m$ by $m-1$, we then obtain (10). \qed

\proclaim{Corollary 4} For any $m,n\in\Z^+$ we have
$$\aligned&\sum_{k=0}^m\bi mkE_{m-k}E_{n-1+k}
\\=&\f{2^{m+1}}n\sum_{k=1}^n\bi nk(2^n-2^{k+1})B_{n-k}\f{(1-2^{m+k+1})B_{m+k+1}}{m+k+1}.
\endaligned\tag11$$
\endproclaim
\Proof. Take $x=y=1/2$ and $z=0$ in (6), and note that $B_n(1/2)=(2^{1-n}-1)B_n$
(cf. [S1]) and also $(-1)^nB_n=B_n$ unless $n=1$. \qed

Theorem 1 will be proved in the next section.
In the proof we use the technique developed in [PS]
together with Theorem 0.

\heading 2. Proof of Theorem 1\endheading

To prove Theorem 1, we need two lemmas.

\proclaim{Lemma 1 {\rm ([PS, Lemma 2.1])}} Let $P(x),Q(x)\in \C[x]$ where $\C$ is the field of complex numbers.

{\rm (i)} We have
$$\Delta(P(x)Q(x))=P(x+1)\Delta(Q(x))+\Delta(P(x))Q(x)$$
and
$$\align\Delta^*(P(x)Q(x))=&\Delta(P(x))Q(x+1)+P(x)\Delta^*(Q(x))
\\=&P(x+1)\Delta^*(Q(x))-\Delta(P(x))Q(x).
\endalign$$

{\rm (ii)} If $\Delta(P(x))=\Delta(Q(x))$, then $P'(x)=Q'(x)$.
 If $\Delta^*(P(x))=\Delta^*(Q(x))$, then $P(x)=Q(x)$.
\endproclaim

\proclaim{Lemma 2} Let $a_0,a_1,\ldots$ be a sequence of complex
numbers, and set
$$A_l(t)=\sum_{k=0}^l\bi lk(-1)^ka_kt^{l-k}$$ for
$l=0,1,2,\ldots$. Then, for any $m,n\in\N$, we have
$$\align&\sum_{k=0}^n\bi nk\f{x^{m+k+1}}{m+k+1}A_{n-k}(y)+(-1)^m\f{A_{m+n+1}(y)}{(m+n+1)\bi{m+n}n}
\\&\quad\ =\sum_{k=0}^m\f{\bi mk}{\bi{n+k}k}(-1)^kx^{m-k}\f{A_{n+k+1}(x+y)}{n+k+1}.
\endalign$$
\endproclaim
\Proof. By [S2], $A_{n+1}'(t)=(n+1)A_n(t)$ and
$$\sum_{k=0}^n\bi nkA_{n-k}(y)z^k=A_n(y+z).$$
Therefore
$$\align &\sum_{k=0}^n\bi nkA_{n-k}(y)\f{x^{m+k+1}}{m+k+1}
=\sum_{k=0}^n\bi nkA_{n-k}(y)\int_0^xt^{m+k}\t dt
\\=&\int_0^xt^m\sum_{k=0}^n\bi nkA_{n-k}(y)t^k\t dt
=\int_0^xt^mA_n(y+t)\t dt
\\=&t^m\f{A_{n+1}(y+t)}{n+1}\bigg|_{t=0}^x-\f m{n+1}\int_0^xt^{m-1}A_{n+1}(y+t)\t dt
\\=&t^m\f{A_{n+1}(y+t)}{n+1}\bigg|_{t=0}^x-\f m{n+1}\cdot\f{t^{m-1}A_{n+2}(y+t)}{n+2}\bigg|_{t=0}^x
\\&+\f m{n+1}\cdot\f{m-1}{n+2}\int_0^xt^{m-2}A_{n+2}(y+t)\t dt
\\=&\cdots=\sum_{k=0}^m(-1)^k\f{m(m-1)\cdots(m-k+1)}{(n+1)\cdots(n+k+1)}t^{m-k}A_{n+k+1}(y+t)\bigg|_{t=0}^x
\\=&\sum_{k=0}^m(-1)^k\f{\bi mk}{\bi {n+k}k}x^{m-k}\f{A_{n+k+1}(x+y)}{n+k+1}
-(-1)^m\f{\bi mm}{\bi{n+m}m}\cdot\f{A_{m+n+1}(y)}{m+n+1}.
\endalign$$
This proves the desired identity. \qed

\medskip
\noindent {\tt Proof of Theorem 1}. We fix $y$ and view $z=1-x-y$ as a function in $x$.

Let $P_{m,n}(x)$
denote the left hand side of (1). Then, with the help of Lemma 1, $\Delta(P_{m,n}(x))$ coincides with
$$\align &(-1)^m\sum_{k=0}^m\bi mk\Delta\(\f{B_{m-k+1}(x)}{m-k+1}\)\f{B_{n+k+1}(y)}{n+k+1}
\\&+(-1)^n\sum_{k=0}^n\bi nk\Delta\(\f{B_{n-k+1}(x)}{n-k+1}\cdot\f{B_{m+k+1}(z)}{m+k+1}\)
\\=&(-1)^m\sum_{k=0}^m\bi mkx^{m-k}\f{B_{n+k+1}(y)}{n+k+1}
+(-1)^n\sum_{k=0}^n\bi nkx^{n-k}\f{B_{m+k+1}(z)}{m+k+1}
\\&+(-1)^n\sum_{k=0}^n\bi nk\f{B_{n-k+1}(x+1)}{n-k+1}\cdot\f{B_{m+k+1}(z-1)-B_{m+k+1}(z)}{m+k+1}.
\endalign$$
In view of Theorem 0 and the above,
$$\align&\Delta(P_{m,n}(x))-\f{(-x)^{m+n+1}}{(m+n+1)\bi {m+n}n}
\\=&\f{(-1)^{n+1}}{n+1}\sum_{k=0}^n\bi{n+1}kB_{n+1-k}(x+1)(z-1)^{m+k}
\\=&\f{(-1)^{n+1}}{n+1}(z-1)^m\l(B_{n+1}(x+1+z-1)-(z-1)^{n+1}\r)
\\=&\f{(z-1)^m}{n+1}\l((-1)^{n+1}B_{n+1}(1-y)-(1-z)^{n+1}\r)
\\=&(-1)^m(x+y)^m\f{B_{n+1}(y)}{n+1}-\f{(-1)^m}{n+1}(x+y)^{m+n+1}.
\endalign$$
Therefore $\Delta(P_{m,n}(x))=\Delta(Q_{m,n}(x))$, where
$$\align Q_{m,n}(x)=&\f{(-1)^{m+n+1}}{(m+n+1)\bi {m+n}n}\cdot\f{B_{m+n+2}(x)}{m+n+2}
\\&+(-1)^m\f{B_{m+1}(x+y)}{m+1}\cdot\f{B_{n+1}(y)}{n+1}-\f{(-1)^m}{n+1}\cdot\f{B_{m+n+2}(x+y)}{m+n+2}.
\endalign$$
By Lemma 1(ii) we must have $P'_{m,n}(x)=Q'_{m,n}(x)$.

It is easy to see that
$$\align Q'_{m,n}(x)=&\f{(-1)^{m+n+1}}{(m+n+1)\bi {m+n}n}B_{m+n+1}(x)
\\&+(-1)^mB_m(x+y)\f{B_{n+1}(y)}{n+1}-\f{(-1)^m}{n+1}B_{m+n+1}(x+y)
\\=&\f{(-1)^{m+n+1}B_{m+n+1}(x)}{(m+n+1)\bi {m+n}n}
+B_m(z)\f{B_{n+1}(y)}{n+1}+\f{(-1)^n}{n+1}B_{m+n+1}(z).
\endalign$$
Also,
$$\align &(-1)^n\(P_{m,n}'(x)-(-1)^m\sum_{k=0}^m\bi mkB_{m-k}(x)\f{B_{n+k+1}(y)}{n+k+1}\)
\\=&\sum_{k=0}^n\bi nk\(B_{n-k}(x)\f{B_{m+k+1}(z)}{m+k+1}-\f{B_{n-k+1}(x)}{n-k+1}B_{m+k}(z)\)
\\=&\sum_{k=1}^n\bi n{k-1}B_{n-k+1}(x)\f{B_{m+k}(z)}{m+k}+\f{B_{m+n+1}(z)}{m+n+1}
\\&-\f{B_{n+1}(x)}{n+1}B_m(z)-\sum_{k=1}^n\bi nk\f{m+k}{n-k+1}B_{n-k+1}(x)\f{B_{m+k}(z)}{m+k}
\\=&\f{B_{m+n+1}(z)}{m+n+1}-\f{B_{n+1}(x)}{n+1}B_m(z)
-m\sum_{k=1}^n\bi nk\f{B_{n-k+1}(x)}{n-k+1}\cdot\f{B_{m+k}(z)}{m+k},
\endalign$$
where in the last step we note that $\bi nk\f k{n-k+1}=\bi n{k-1}$ for $k=1,2,\ldots,n$.
Observe that $P_{m+1,n}'(x)$ coincides with
$$\align &(-1)^{m+1}\sum_{k=0}^{m+1}\bi {m+1}kB_{m-k+1}(x)\f{B_{n+k+1}(y)}{n+k+1}
+(-1)^n\f{B_{m+n+2}(z)}{m+n+2}
\\&-(-1)^n(m+1)\sum_{k=0}^n\bi nk\f{B_{n-k+1}(x)}{n-k+1}\cdot\f{B_{m+k+1}(z)}{m+k+1}
\\=&-(m+1)P_{m,n}(x)+(-1)^{m+1}\f{B_{m+n+2}(y)}{m+n+2}+(-1)^n\f{B_{m+n+2}(z)}{m+n+2},
\endalign$$
On the other hand, $Q_{m+1,n}'(x)$ equals
$$\f{(-1)^{m+n}B_{m+n+2}(x)}{(m+n+2)\bi{m+n+1}n}
+\f{B_{m+1}(z)B_{n+1}(y)+(-1)^nB_{m+n+2}(z)}{n+1}.$$
Now it is clear that the equality $P_{m+1,n}'(x)=Q_{m+1,n}'(x)$ yields (1).

Next we come to prove (2). Let $L(x)$ denote the left hand side of (2). Then
$$\align \Delta^*(L(x))=&(-1)^m\sum_{k=0}^m\bi mk\Delta^*(E_{m-k}(x))\f{B_{n+k+1}(y)}{n+k+1}
\\&+(-1)^n\sum_{k=0}^n\bi nk\Delta^*\(E_{n-k}(x)\f{B_{m+k+1}(z)}{m+k+1}\).
\endalign$$
Applying Lemma 1 and the first identity in Theorem 0, we then obtain that
$$\align \Delta^*(L(x))=&(-1)^m\sum_{k=0}^m\bi mk2x^{m-k}\f{B_{n+k+1}(y)}{n+k+1}
\\&+(-1)^n\sum_{k=0}^n\bi nk2x^{n-k}\f{B_{m+k+1}(z)}{m+k+1}
\\&+(-1)^n\sum_{k=0}^n\bi nkE_{n-k}(x+1)\f{B_{m+k+1}(z-1)-B_{m+k+1}(z)}{m+k+1}
\\=&\f{2(-x)^{m+n+1}}{(m+n+1)\bi {m+n}n}-(-1)^n\sum_{k=0}^n\bi nkE_{n-k}(x+1)(z-1)^{m+k}
\\=&\f{2(-x)^{m+n+1}}{(m+n+1)\bi {m+n}n}-(-1)^n(z-1)^mE_n(x+1+z-1)
\\=&\f{2(-x)^{m+n+1}}{(m+n+1)\bi {m+n}n}-(-1)^m(x+y)^mE_n(y).
\endalign$$
Therefore
$$\Delta^*(L(x))=\Delta^*\(\f{(-1)^{m+n+1}E_{m+n+1}(x)}{(m+n+1)\bi{m+n}n}
-\f{(-1)^m}2E_m(x+y)E_n(y)\).$$
In view of Lemma 1(ii), we have
$$L(x)=\f{(-1)^{m+n+1}E_{m+n+1}(x)}{(m+n+1)\bi{m+n}n}
-\f{(-1)^m}2E_m(x+y)E_n(y)$$
which is equivalent to the desired (2).

Finally let us turn to prove (3). Let
$f(x)$ denote the left hand side of (3).
By Lemma 1,
$$\align\Delta^*(f(x))=&\f{(-1)^m}2\sum_{k=0}^m\bi mk\Delta^*(E_{m-k}(x))\f{E_{n+k+1}(y)}{n+k+1}
\\&-(-1)^{n}\sum_{k=0}^n\bi nk\Delta^*\(\f{B_{n-k+1}(x)}{n-k+1}\cdot\f{E_{m+k+1}(z)}{m+k+1}\)
\\=&(-1)^m\sum_{k=0}^m\bi mkx^{m-k}\f{E_{n+k+1}(y)}{n+k+1}
\\&+(-1)^{n}\sum_{k=0}^n\bi nkx^{n-k}\f{E_{m+k+1}(z)}{m+k+1}
\\&-(-1)^{n}\sum_{k=0}^n\bi nk\f{B_{n-k+1}(x+1)}{n-k+1}\cdot\f{2(z-1)^{m+k+1}}{m+k+1}.
\endalign$$
In view of the first identity
in Theorem 0 with Bernoulli polynomials replaced by corresponding Euler polynomials, we have
$$\Delta^*(f(x))=\f{(-x)^{m+n+1}}{(m+n+1)\bi{m+n}n}+\f{(-1)^{n+1}}{n+1}2R$$
where
$$\align R=&\sum_{k=0}^n\bi{n+1}kB_{n+1-k}(x+1)\f{(z-1)^{m+k+1}}{m+k+1}
\\=&-\f{(z-1)^{m+n+2}}{m+n+2}-(-1)^m\f{B_{m+n+2}(x+1)}{(m+n+2)\bi{m+n+1}{n+1}}
\\&+\sum_{k=0}^m\f{\bi mk}{\bi{n+1+k}k}(-1)^k(z-1)^{m-k}\f{B_{n+k+2}(x+1+z-1)}{n+k+2}
\endalign$$
by applying Lemma 2 with $a_k=(-1)^kB_k$. Therefore
$$\align \Delta^*(f(x))=&-\sum_{k=0}^m\f{\bi mk}{\bi{n+k+1}k}\cdot\f{2(z-1)^{m-k}}{n+1}\cdot\f{B_{n+k+2}(y)}{n+k+2}
\\&+\f{(-1)^{m+n}}{(m+n+1)\bi{m+n}n}\(\f{2B_{m+n+2}(x+1)}{m+n+2}-x^{m+n+1}\)
\\&+\f{(-1)^n}{n+1}\cdot\f{2(z-1)^{m+n+2}}{m+n+2}.
\endalign$$
Let $g(x)$ denote the right hand side of (3). Clearly $\Delta^*(g(x))$ also
coincides with the right hand side of the last equality. Thus
$\Delta^*(f(x))=\Delta^*(g(x))$ and hence $f(x)=g(x)$ as desired. We are done. \qed

\enddocument